\documentclass[12pt]{amsart}
\usepackage{amssymb}
\usepackage{hyperref}
\usepackage{amsmath}
\usepackage{amscd}
\usepackage{graphicx}
\usepackage{xcolor}
\usepackage{bm}
\usepackage{comment}

\DeclareMathSymbol{\shortminus}{\mathbin}{AMSa}{"39}

\newcommand{\group}[1]{\mathrm{#1}}

\newcommand{\rep}[1]{\rep{#1}}

\newcommand{\Aut}{\operatorname{Aut}}

\newcommand{\Tr}{\operatorname{Tr}}

\newcommand{\C}{\mathbb{C}}
\newcommand{\D}{\mathbb{D}}

\newcommand{\N}{\mathbb{N}}

\newcommand{\R}{\mathbb{R}}

\newcommand{\V}{\mathbf{V}}
\newcommand{\W}{\mathbf{W}}

\newcommand{\Y}{\mathsf{Y}}

\newtheorem{thm}{Theorem}[section]
\newtheorem{prop}[thm]{Proposition}

\newtheorem{cor}[thm]{Corollary}

\theoremstyle{definition}

\title[Quasimodular Asymptotics]{Quasimodular Asymptotics of Spherical Integrals}

\author{Jonathan Novak}

\begin{document}

\begin{abstract}
We show that the spherical integral of the Circular Unitary Ensemble converges 
in expectation to Euler's generating function for integer partitions, 
and that subleading corrections to this high-dimensional limit are quasimodular forms. 
\end{abstract}

\maketitle

\section{Introduction}
Let $\group{U}_N=\Aut \C^N$ be the unitary group with its normalized Haar measure, $\mathrm{d}U$.
In this paper we take a new point of view on the spherical integral

    \begin{equation}
    \label{eqn:SphericalIntegral}
        K_N(q,A,B) = \int e^{qN\Tr AUBU^*} \mathrm{d}U
    \end{equation}

\noindent
of Harish-Chandra \cite{HC} and Itzykson-Zuber \cite{IZ} by considering it as an
observable of the Circular Unitary Ensemble \cite{Mehta}. More precisely, 
taking $A=A_N$ and $B=A_N^*$ with $A_N$ Haar-distributed in $\group{U}_N$, 
the determinantal evaluation \cite{HC,IZ} of \eqref{eqn:SphericalIntegral} becomes
 
    \begin{equation}
    \label{eqn:RandomHCIZ}
       K_N(q,A_N,A_N^*) = \frac{1!2! \dots (N-1)!}{(qN)^{N \choose 2}}
        \frac{\det [e^{qNa_i\bar{a}_j}]}{\prod_{i<j} |a_i-a_j|^2}
    \end{equation}

\noindent
with $a_1,\dots,a_N$ the random eigenvalues of $A_N$. Like the characteristic polynomial of $A_N$,
this random entire function of $q \in \C$ is visibly sensitive to eigenvalue repulsion. We show here
that for $|q|< e^{-1}$ the expected spherical integral of the CUE converges without normalization
to Euler's generating function for integer partitions:

    \begin{equation}
    \label{eqn:EulerLimit}
       \lim_{N \to \infty}\iint e^{qN\Tr AUA^*U^*} \mathrm{d}U\mathrm{d}A  = \prod_{n=1}^\infty \frac{1}{1-q^n}.
    \end{equation}

\noindent
Thus for example we have 

    \begin{equation}
        \lim_{N \to \infty}\iint e^{e^{-\pi}N\Tr AUA^*U^*} \mathrm{d}U\mathrm{d}A = 
        \frac{2^{7/8}\pi^{3/4}}{e^{\pi/24}\Gamma(\frac{1}{4})},
    \end{equation}

\noindent
where the evaluation of $\prod_{n=1}^\infty (1-e^{-n\pi})$ is due to Ramanujan \cite{Berndt}. 
In fact, we prove a much more precise result: the principal logarithm $L_N(q)$ of the expected spherical integral of $A_N$
admits a complete $N \to \infty$ asymptotic expansion whose subleading 
terms are quasimodular forms. This gives a new point of contact between 
random matrices and number theory. Let $\|F\|_\varepsilon$ denote the 
maximum modulus of a continuous function $F(q)$ on the closure of 
$\D(\varepsilon) = \{q \in \C \colon |q|<\varepsilon\}$, and let $\D=\D(1)$.

    \begin{thm}
    \label{thm:Main}
   For any $\varepsilon \in [0,e^{-1})$ and any $m \in \N$ we have

            \begin{equation*}
                \lim_{N \to \infty} N^{2m-2} \left\|L_N - \sum_{g=1}^m N^{2-2g} F_g\right\|_\varepsilon =0,
            \end{equation*}

        \noindent
    where $F_g(q) = \sum_{d=1}^\infty \frac{q^d}{d!} F_g^d$
    is an analytic function on $\D$ and $F_g^d$ counts configurations
    $(\pi_1,\pi_2,\tau_1,\dots,\tau_{2g-2})$ of permutations from
    $\group{S}^d=\Aut\{1,\dots,d\}$ such that:

            \begin{enumerate}

                \smallskip
                \item 
                The product $\pi_1\pi_2\pi_1^{-1}\pi_2^{-1}\tau_1\dots\tau_{2g-2}$ equals the identity
                permutation;

                \smallskip
                \item 
                The permutations $\tau_1,\dots,\tau_{2g-2}$ are transpositions;

                 \smallskip
                \item 
                The group generated by $\pi_1,\pi_2,\tau_1,\dots,\tau_{2g-2}$
                acts transitively on $\{1,\dots,d\}$;

                 \smallskip
                \item 
                Writing $\tau_i = (s_i\ t_i)$ with $1 \leq s_i < t_i \leq d$ we have
                $t_1 \leq \dots \leq t_{2g-2}$.
            \end{enumerate}

    \end{thm}

It is not difficult to see that $F_1(q)$ is indeed the logarithm of Euler's generating function,

    \begin{equation}
        F_1(q) = \sum_{n=1}^\infty \frac{q^n}{1-q^n} n^{-1}.
    \end{equation}

\noindent
Theorem \ref{thm:Main} makes quasimodularity of 
subleading corrections to the limit $L_N \to F_1$ quite believable. Let $C_g^d$ be the number of 
configurations as in Theorem \ref{thm:Main} but subject only to the 
first three conditions therein, and note that $C_1^d=F_1^d.$ For $g \geq 2,$ 
it is a theorem of Dijkgraaf \cite{Dijkgraaf} and Kaneko-Zagier \cite{KanZag} that
the generating function $C_g(q)=\sum_{d=1}^\infty \frac{q^d}{d!}C_g^d$ is a polynomial in the Eisenstein series

    \begin{equation}
        \label{eqn:Eisenstein}
        E_{2k}(q) = 1 + \frac{2}{\zeta(1-2k)}\sum_{n=1}^\infty \frac{q^n}{1-q^n} n^{2k-1}, \quad k \in \N,
    \end{equation}

\noindent
which with the exception of $E_2$ are Fourier expansions of modular forms. 
More recently, the quasimodularity of a large class of generating functions
containing both $C_g(q)$ and $F_g(q)$ was established by Hahn, van Ittersum, and 
Leid \cite{HIL} using the Bloch-Okounkov theorem \cite{BlOk,Zagier}.
Since the algebra $\C[E_2,E_4,E_6,\dots]$ of quasimodular forms
is freely generated by $E_2,E_4,E_6$, Theorem \ref{thm:Main} and quasimodularity of $F_g(q)$ for $g \geq 2$ together
imply that the subleading asymptotics of $L_N(q)$ 
can be expressed solely in terms of these three Eisenstein series. For example,

    \begin{equation}
       L_N = F_1 + \frac{5E_2^3-3E_2E_4-2E_6+45E_2^2+18E_4+90E_2-153}{51840N^2} + O\left(\frac{1}{N^4} \right)
    \end{equation}

\noindent
as $N \to \infty$, uniformly on compact subsets of $\D(e^{-1})$.

Theorem \ref{thm:Main} has a dual interpretation in terms of the statistical mechanics
of partitions. Let $\Y_N^d$ denote the set of Young diagrams
with exactly $d$ cells and at most $N$ rows. For each $\lambda \in \Y_N^d$, let 
$\V^\lambda$ and $\W_N^\lambda$ denote the corresponding irreducible representations
of $\group{S}^d$ and $\group{U}_N$, respectively, and let
$s_\lambda(A)$ denote the character of $A \in \group{U}_N$ acting in $\W_N^\lambda$.
From the character expansion \cite[Equation (3.20)]{IZ} of $K_N(q,A,B)$, we obtain the 
Maclaurin expansion of the random entire function $K_N(q,A_N,A_N^*)$ as

                \begin{equation}
                \label{eqn:CharacterHCIZ}
                    K_N(q,A_N,A_N^*) = 1 + \sum_{d=1}^\infty \frac{q^d}{d!}N^d
                    \sum_{\lambda \in \Y_N^d} |s_\lambda(A_N)|^2 \frac{\dim \V^\lambda}{\dim \W_N^\lambda}.
                \end{equation}

\noindent
Thus, $K_N(q,A_N,A_N^*)$ takes positive values for positive arguments and
can be viewed as the partition function of a random 
probability measure on the set $\Y_N$ of Young diagrams with at most $N$ rows.
In the event $A_N$ is a scalar multiple of the identity, this random measure becomes
the nonrandom Poissonized Schur-Weyl measure on $\Y_N$, whose $N \to \infty$ asymptotics have been 
analyzed by Borodin and Olshanski \cite{BoOl}. 
From this point of view, Theorem \ref{thm:Main} is an $N \to \infty$
asymptotic expansion of the annealed free energy 

      \begin{equation}
        L_N(q) =  \log \int K_N(q,A,A^*) \mathrm{d}A
    \end{equation}

\noindent
in the Poissonized Schur-Weyl model with CUE disorder. 

The spherical integral $K_N(q,A,B)$ belongs to a large class of multivariate 
special functions: hypergeometric functions of two complex matrix arguments $A$ and $B$. 
These functions were introduced by James \cite{James}
in multivariate statistics and have been much studied in connection with spiked Hermitian matrix models \cite{JO,PMC}.
Recently, the present author proposed to view arbitrary
hypergeometric functions of $A=A_N$ and $B=A_N^*$ as a highly structured
but completely unexplored family of spectral observables of the CUE, and showed
that a general version of Theorem \ref{thm:Main} holds modulo  
a nontrivial interchange of limits \cite{Novak:hyper}.
The detailed estimates presented here rigorously justify this 
limit swap for the expected spherical integral of the CUE.

It is interesting to compare the asymptotics of the
integral $K_N(q)=\int K_N(q,A,A^*)\mathrm{d}A$ with those of the integrand. 
The logarithm of the integrand admits an $N \to \infty$ 
asymptotic expansion of the form

    \begin{equation}
    \label{eqn:NoAverage}
        \log K_N(q,A,A^*) \sim \sum_{g=0}^\infty N^{2-2g}F_g(q,A,A^*)
    \end{equation}

\noindent
which holds uniformly over $A \in \group{U}_N$ provided $|q|$ is 
sufficiently small, with the coefficients $F_g(q,A,A^*)$ 
being generating functions for transposition factorizations similar to those appearing
in Theorem \ref{thm:Main} and subject to the same monotonicity constraint
on transposition factors \cite{GGN3,Novak:Fourier}. This constraint defines monotone Hurwitz
theory \cite{GGN1,GGN2}, a combinatorial desymmetrization of the classical branch of 
enumerative geometry concerned with counting branched covers of Riemann surfaces with 
given branch locus and ramification profile.
The coefficients $F_g(q,A,A^*)$ in the asymptotic expansion \eqref{eqn:NoAverage} are generating
functions for monotone double Hurwitz numbers over the Riemann sphere, whereas the coefficients $F_g(q)$
in Theoem \ref{thm:Main} are generating functions for monotone simple Hurwitz numbers over
an elliptic curve. Thus, the drop in leading order from $O(N^2)$ to
$O(1)$ after averaging $A$ over $\group{U}_N$ is a topological manifestation of eigenvalue 
repulsion in the CUE. The proof of Theorem \ref{thm:Main} presented here is a standalone 
argument which does not rely on any results from the existing literature on high-dimensional asymptotics
of spherical integrals \cite{GuiMai,GuiNov,GuiZei} --- in fact, it is a helpful guide to the more intricate 
arguments in \cite{Novak:Fourier}.

\section{First Order Asymptotics}
\label{sec:FirstOrder}
In this section we prove Theorem \ref{thm:Main} for $m=1$.
    
    \subsection{Convergence of derivatives}
    Let us write 

        \begin{equation}
            K_N(q) = \iint e^{qN\Tr AUA^*U^*} \mathrm{d}U\mathrm{d}A= 1 + \sum_{d=1}^\infty \frac{q^d}{d!}K_N^d
        \end{equation}

    \noindent
    for the expected spherical integral of the CUE, 
    viewed as an entire function of $q \in \C$. As an immediate consequence of 
    \eqref{eqn:CharacterHCIZ} and $\int |s_\lambda(A)|^2 \mathrm{d}A=1$, we have

         \begin{equation}
        \label{eqn:YoungSum}
            K_N^d = N^d\sum_{\lambda \in \Y_N^d} \frac{\dim \V^\lambda}{\dim \W_N^\lambda}.
        \end{equation}

    \noindent
    This formula can be used to establish the $N \to \infty$ convergence of 
    $K_N^d.$
    
    For any Young diagram $\lambda$, the content $c(\Box)$ of a given cell $\Box \in \lambda$
    is its column index minus its row index: filling the cells of $\lambda$ with their 
    contents produces a $\lambda$-shaped corner of the infinite Toeplitz matrix $[j-i]_{i,j=1}^\infty$.
    The content polynomial \cite[Chapter 1, Section 1, Example 11]{Macdonald} of $\lambda$
    is the product

        \begin{equation}
        \label{eqn:ContentPolynomial}
            x^\lambda = \prod_{\Box \in \lambda} (x+c(\Box)).
        \end{equation}

    \noindent
    As $\lambda$ ranges over the set $\Y^d$ of all Young diagrams with exactly $d$ cells, 
    $x^\lambda$ interpolates between the falling factorial $x^{\underline{d}}=x(x-1) \dots (x-d+1)$,
    which corresponds to a column of $d$ cells, and the rising factorial 
    $x^{\overline{d}}=x(x+1) \dots (x+d-1)$, which corresponds to a row of $d$ cells. 
    It is clear that if $x > d-1$ then $x^\lambda >0$ for all $\lambda \in \Y^d$, 
    and moreover we have $x^\lambda \leq x^\mu$ if $\lambda \leq \mu$
    with respect to dominance order on $\Y^d$.
    From the hook-length and hook-content formulas \cite[Chapter 1, Section 3, Example 4]{Macdonald} for the 
    dimensions of $\V^\lambda$ and $\W_N^\lambda$, we have

        \begin{equation}
        \frac{\dim \V^\lambda}{\dim \W_N^\lambda} = \frac{d!}{N^\lambda}, \quad \lambda \in \Y_N^d.
     \end{equation}

    \noindent
    Thus, \eqref{eqn:YoungSum} may be restated as follows.

    \begin{prop}
        \label{prop:YoungSum}
        For each $d \in \N$, we have

            \begin{equation*}
                K_N^d = d! \sum_{\lambda \in \Y_N^d} \frac{N^d}{N^\lambda} = 
                d! \sum_{\lambda \in \Y_N^d} \prod_{\Box \in \lambda} \frac{1}{1+\frac{c(\Box)}{N}}.
            \end{equation*}
    \end{prop}

    From Proposition \ref{prop:YoungSum}, we immediately have the convergence

        \begin{equation}
        \label{eqn:DerivativeConvergence}
            \lim_{N \to \infty} K_N^d =d!|\Y^d|,
        \end{equation}

    \noindent
    as first obtained in \cite{Novaes}. Let 
    $H_1^d= d!|\Y^d|$, where the subscript is present because
    in Section \ref{sec:Stable} we will obtain subleading 
    corrections $H_2^d,H_3^d,\dots$ to the above limit.

    \subsection{Uniform boundedness}
    We now use Proposition \ref{prop:YoungSum} to show that $K_N(q)$ is 
    uniformly bounded in a neighborhood of $q=0.$ We prove this by partitioning
    the Maclaurin expansion of $K_N(q)$ into a stable part and an unstable part,
    and bounding the modulus of each.

        \begin{prop}
            \label{prop:Uniform}
            For any $\gamma \in [0,e^{-1}),$ we have $\sup_{N \in \N} \|K_N\|_\gamma < \infty.$
        \end{prop}

        \begin{proof}
        By the stable part of $K_N$, we mean its degree $N$ Maclaurin 
        polynomial 

            \begin{equation*}
                K_{NN}(q) = 1 + \sum_{d=1}^N \frac{q^d}{d!}K_N^d.
            \end{equation*}

         \noindent
    The derivatives $K_N^1,\dots,K_N^N$ are stable in the sense
    that, for $1 \leq d \leq N$, Proposition \ref{prop:YoungSum} 
    gives $K_N^d$ as a sum over the full set $\Y^d$ of Young diagrams with $d$ cells. 
        For any $1 \leq d \leq N$, the function on $\Y^d$ defined by 
        $\lambda \mapsto N^\lambda$ is minimized by taking $\lambda$
        to be a column of $d$ cells and the minimum value is the falling
        factorial $N^{\underline{d}} = N(N-1) \dots (N-d+1),$
        giving the bound 

            \begin{equation*}
               K_N^d = d! \sum_{\lambda \in \Y^d}\frac{N^d}{N^\lambda} \leq d!|\Y^d|\frac{N^d}{N^{\underline{d}}}.
            \end{equation*}

        \noindent
        The sequence

            \begin{equation*}
                \frac{N^d}{N^{\underline{d}}} = \frac{N}{N} \frac{N}{N-1} \dots \frac{N}{N-d+1},
                \quad N \geq d,
            \end{equation*}

        \noindent
        is monotonically decreasing in $N$ with initial term $\frac{d^d}{d!} < e^d$,
        so for $\gamma \in [0,e^{-1})$ we have

            \begin{equation*}
               \|K_{NN}\|_\gamma < 1 + \sum_{d=1}^N \gamma^d \sum_{\lambda \in \Y^d}\frac{N^d}{N^\lambda} <
                 1 + \sum_{d=1}^N (\gamma e)^d|\Y_d| < \prod_{n=1}^\infty \frac{1}{1-\gamma^ne^n}.
            \end{equation*}

    By the unstable part of $K_N(q)$, we mean the tail sum

        \begin{equation}
            \label{eqn:StableFirstOrder}
            K_N(q) - K_{NN}(q) =\sum_{d=N+1}^\infty \frac{q^d}{d!}K_N^d.
        \end{equation}

    \noindent
    The higher derivatives $K_N^d$ are unstable in the sense that
    for $d>N$ the sum in Proposition \ref{prop:YoungSum} is over the proper 
    $N$-dependent subset $\Y_N^d$ of $\Y^d$.
    For $d > N$, we have $d=rN+s$ for some $r \in \N$ and $s \in \{0,1,\dots,N-1\}$,
            with $s>0$ if $r=1$.
            The function $\lambda \to N^\lambda$ is minimized by taking $\lambda$ to be the 
            diagram consisting of $r$ columns of length $N$ together with a final column of length $s$,
            and for this diagram

                \begin{equation*}
                \begin{split}
                    N^\lambda &= \frac{N!}{0!} \frac{(N+1)!}{1!} \frac{(N+2)!}{2!} \dots \frac{(N+r-1)!}{(r-1)!}
                    (N+r)^{\underline{s}} \\
                    &= (N!)^r {N + 0 \choose 0} {N+1 \choose 1} \dots {N+r-1 \choose r-1} (N+r)^{\underline{s}}.
                \end{split}
                \end{equation*}

            \noindent
            This gives 
            
                \begin{equation*}
                    N^\lambda \geq (N!)^r (N+r)^{\underline{s}} > (N!)^{\frac{d}{N}}> \frac{N^d}{e^d},
                \end{equation*}

            \noindent
            so for $\gamma \in [0,e^{-1})$ we have

                \begin{equation*}
                    \|K_N - K_{NN}\|_\gamma < \sum_{d=N+1}^\infty (\gamma e)^d |\Y^d|.
                \end{equation*}

            \noindent
            Using the standard upper bound $|\Y^d| \leq e^{\pi\sqrt{\frac{2d}{3}}}$ on partitions,
            a straightforward calculation (see e.g. \cite[Lemma 2.4]{LeMa}) yields

                \begin{equation*}
                    \sum_{d=N+1}^\infty (\gamma e)^d |\Y^d| \leq C(\gamma) e^{\frac{\log (\gamma e)}{2}N},
                \end{equation*}

            \noindent
            with $C(\gamma) \geq 0$ independent of $N$.
        \end{proof}

    \subsection{Uniform convergence}
    We now prove the convergence of the expected spherical integral $K_N(q)=\iint e^{qN\Tr AUA^*U^*}\mathrm{d}U\mathrm{d}A$ 
    to the Euler product $H_1(q)=\prod_{n=1}^\infty (1-q^n)^{-1}.$

        \begin{thm}
            \label{thm:ConvergenceToEuler}
            We have $K_N \to H_1$ as $N \to \infty$, uniformly on compact subsets of $\D(e^{-1}).$
        \end{thm}

        \begin{proof}
            Fix $0 \leq \varepsilon < \gamma < e^{-1}$. We have

                \begin{equation*}
                    \|K_N-H_1\|_\varepsilon \leq \sum_{d=1}^\infty \frac{\varepsilon^d}{d!}|K_N^d-H_1^d| < \infty.
                \end{equation*}

            \noindent
            By Cauchy's inequality \cite[Section 2.5]{Titchmarsh}, 

                \begin{equation*}
                    \frac{1}{d!} |K_N^d-H_1^d| \leq \frac{\|K_N-H_1\|_\gamma}{\gamma^d} \leq \frac{M}{\gamma^d},
                    \quad d \in \N,
                \end{equation*}

            \noindent
            where $M=\sup_N \|K_N\|_\gamma + \|H_1\|_\gamma$ is finite by Proposition \ref{prop:Uniform}.
            Thus for any $r \in \N$ we have

                \begin{equation*}
                     \|K_N-H_1\|_\varepsilon \leq \sum_{d=1}^r \frac{\varepsilon^d}{d!}|K_N^d-H_1^d| + \left(\frac{\varepsilon}{\gamma}\right)^{r+1}
                     \frac{M}{1-\frac{\varepsilon}{\gamma}}.
                \end{equation*}

            \noindent
            Given arbitrarily small $x>0$, take $r$ large enough so that the tail term is smaller than $\frac{x}{2}.$
            By \eqref{eqn:DerivativeConvergence}, the sum of the other $r$ terms is less than $\frac{x}{2}$
            for all $N$ sufficiently large.
        \end{proof}

Now consider the principal logarithm $L_N(q)=\log K_N(q),$ which is defined and analytic on an open
origin-centered disc of positive radius $\kappa_N$ equal to the minimum modulus of a zero of $K_N(q).$ 
For any $\gamma \in (0,e^{-1}),$ Theorem \ref{thm:ConvergenceToEuler} and the Hurwitz theorem on 
zeros of uniform limits \cite[Section 3.45]{Titchmarsh} together imply the existence of $N(\gamma) \in \N$ such that
$K_N(q)$ is nonvanishing on the closure of $\D(\gamma)$ for all $N \geq N(\gamma).$ Thus, considering
$L_N(q)$ as a sequence of analytic functions on $\D(\gamma)$ indexed by $N \geq N(\gamma)$, we have
$L_N \to F_1=\log H_1$ uniformly on compact subsets of $\D(\gamma).$

In order to complete the proof of Theorem \ref{thm:Main} in the case $m=1$ it remains only to show that 
the Macluarin expansion of the logarithme of the Euler product is

    \begin{equation}
        F_1(q) = \log \prod_{n=1}^\infty \frac{1}{1-q^n} = \sum_{d=1}^\infty \frac{q^d}{d!} F_1^d,
    \end{equation}

\noindent
with $F_1^d$ the number of commuting 
pairs of permutations of degree $d$ which generate a transitive subgroup of
$\group{S}^d$. By the Exponential Formula \cite[Chapter 5]{Stanley:EC2}, this is equivalent to 
showing that the Euler product itself is an exponential generating function enumerating commuting pairs of permutations. 
This follows directly from the orbit-stabilizer theorem, which can be used to show that
the number of commuting pairs of elements in any finite group is the order of the group times its class
number \cite{Burnside}.

\section{Stability and Concentration}
\label{sec:Stable}
We have now proved Theorem \ref{thm:Main} in the case $m=1$, 
and in Section \ref{sec:HigherOrder} we will give a complete $N \to \infty$ asymptotic
expansion for the principal logarithm $L_N(q)=\log K_N(q)$. 
In this section, we discuss the $N=\infty$ limit of this approximation problem.

    \subsection{Stable derivatives}
    For any integers $1 \leq d \leq N$, Proposition \ref{prop:YoungSum}
    gives the corresponding derivative $K_N^d$
    of $K_N(q)$ at $q=0$ as the value of the rational function 

        \begin{equation}
            K^d(t) = d!\sum_{\lambda \in \Y^d} \prod_{\Box \in \lambda} \frac{1}{1+tc(\Box)}
        \end{equation}

    \noindent
    at $t=1/N.$ This rational function has two dominant singularities 
    in the complex $t$-plane, namely a pair of simple poles at $\pm 1/(d-1)$,
    so its Maclaurin series has radius of convergence $1/(d-1).$ This series
    is given explicitly by

    \begin{equation}
    \label{eqn:RationalMaclaurin}
        K^d(t) = d! \sum_{r=0}^\infty (-t)^r \sum_{\lambda \in \Y^d} h_r(\lambda),
    \end{equation}

    \noindent
    where $h_r(\lambda)$ denotes the complete 
    homogeneous symmetric polynomial of degree $r$ in $d$ variables,

    \begin{equation}
        \label{eqn:Complete}
        h_r(x_1,\dots,x_d) = \sum_{1 \leq i_1 \leq \dots \leq i_r \leq d}
        x_{i_1} \dots x_{i_r},
    \end{equation}

    \noindent
    evaluated on the multiset of contents of the Young diagram $\lambda$. Since $K^d(t)=K^d(-t)$,
    the Maclaurin expansion of $K^d(t)$ is a power series in $t^2$ and we write

    \begin{equation}
        K^d(t) = \sum_{g=1}^\infty t^{2g-2} H_g^d
    \end{equation}

    \noindent
    where 

    \begin{equation}
    \label{eqn:CompleteContents}
        H_g^d = d!\sum_{\lambda \in \Y^d} h_{2g-2}(\lambda).
    \end{equation}

    \noindent
   We thus have the following.

    \begin{thm}
         \label{thm:LeadingDerivatives}
        For any integers $1 \leq d \leq N$ we have

            \begin{equation*}
                K_N^d = \sum_{g=1}^\infty N^{2-2g} H_g^d,
            \end{equation*}

        \noindent
        where the series is absolutely convergent.
    \end{thm}

    \noindent
   Theorem \ref{thm:LeadingDerivatives} gives an $N \to \infty$ asymptotic
    expansion for each fixed derivative $K_N^d$ of $K_N(q)$ at $q=0$, providing all
    subleading corrections to the large $N$ limit $K_N^d \to H_1^d=d!|\Y^d|$.
    In fact, all the coefficients $H_g^d$ are positive integers.
    The coefficients $H_g^d$ in Theorem \ref{thm:LeadingDerivatives} are positive integers.
    Integrality is immediate from \eqref{eqn:CompleteContents}, and we will give two proofs
    of positivity. The first uses a theorem of Hunter on complete symmetric
    polynomials of even degree \cite{Hunter,Tao} and gives bounds on $H_g^d$, while the second is 
    representation-theoretic and gives a combinatorial interpretation of $H_g^d$.

    \begin{thm}[Hunter]
    \label{thm:Hunter}
        For all $g,d \in \N$, we have

            \begin{equation*}
                h_{2g-2}(x_1,\dots,x_d) \geq \frac{(x_1^2+\dots + x_d^2)^{g-1}}{2^{g-1}(g-1)!}
            \end{equation*}

        \noindent
        for every $(x_1,\dots,x_d) \in \R^d$. Moreover, $h_{2g-2}$ is a Schur convex
        function on $\R^d$.
    \end{thm}

Hunter's theorem immediately implies the positivity of $H_g^d,$ and also yields bounds on

    \begin{equation}
        \frac{H_g^d}{H_1^d} = \frac{1}{|\Y^d|} \sum_{\lambda \in \Y^d} h_{2g-2}(\lambda),
    \end{equation}

\noindent
the expected value of $h_{2g-2}(\lambda)$ when $\lambda$ is uniformly distributed
in $\Y^d$. Following \cite[Section 6.1]{GraKnuPat}, we denote the Stirling number
of the second kind by ${m \brace n}.$ 

    \begin{thm}
        \label{thm:StirlingBound}
        For all $g,d \in \N$, we have

            \begin{equation*}
                \frac{d^{g-1}}{2^{g-1}(g-1)!} \leq \frac{H_g^d}{H_1^d} \leq {d-1 + 2g-2 \brace d-1}.
            \end{equation*}
    \end{thm}

    \begin{proof}
        The lower bound is an immediate consequence of the lower bound in Theorem \ref{thm:Hunter}.
       The upper bound follows from the Schur convexity of $h_{2g-2}$, which implies that
        for any $\lambda,\mu \in \Y^d$ with $\mu \leq \lambda$ in dominance order we have
        $h_{2g-2}(\mu) \leq h_{2g-2}(\lambda)$, and therefore 

            \begin{equation*}
                h_{2g-2}(\lambda) \leq h_{2g-2}(1,2,\dots,d-1), \quad \lambda \in \Y^d.
            \end{equation*}

        \noindent
        It is well known \cite[Chapter 1, Section 2, Example 11]{Macdonald} that 
        $h_m(1,2,\dots,n) = {n +m \brace n}$.
    \end{proof}

We will soon use Theorem \ref{thm:StirlingBound} in conjunction with an elementary 
upper bound on Stirling numbers which we now state.

    \begin{prop}
        \label{prop:StirlingBoung} We have ${m \brace n} \leq \frac{1}{n!}n^m< e^nn^{m-n}.$
    \end{prop}

    \begin{proof}
        Since ${m \brace n}$ is the number of partitions of $\{1,\dots,m\}$ into $n$
        disjoint nonempty subsets, we have

            \begin{equation*}
               {m \brace n} = \frac{1}{n!}|\mathrm{Sur}(m,n)|,
            \end{equation*}

        \noindent
        where $\mathrm{Sur}(m,n)$ is the set of surjections $\{1,\dots,m\} \to \{1,\dots,n\}.$
        The claimed bounds thus follow from $|\mathrm{Sur}(m,n)|\leq n^m$ and $n!>e^{-n}n^n.$
    \end{proof}

Our second proof of the positivity of $H_g^d$ is representation-theoretic
and relies on spectral properties of Jucys-Murphy elements in the group algebra
of the symmetric group \cite{Cec,DiGr}.

\begin{thm}
    \label{thm:Monotone}
    The number $H_g^d$ counts configurations
    $(\pi_1,\pi_2,\tau_1,\dots,\tau_{2g-2})$ of permutations from 
    $\group{S}^d=\Aut\{1,\dots,d\}$ such that:

            \begin{enumerate}

                \smallskip
                \item 
                The product $\pi_1\pi_2\pi_1^{-1}\pi_2^{-1}\tau_1\dots\tau_{2g-2}$ equals the identity
                permutation $\iota \in \group{S}^d$;

                \smallskip
                \item 
                The permutations $\tau_1,\dots,\tau_{2g-2}$ are transpositions;


                 \smallskip
                \item 
                Writing $\tau_i = (s_i\ t_i)$ with $1 \leq s_i < t_i \leq d$ we have
                $t_1 \leq \dots \leq t_{2g-2}$.
            \end{enumerate}

    \end{thm}

    \begin{proof}
        Writing 

            \begin{equation*}
                H_g^d = \sum_{\lambda \in \Y^d} \frac{(\dim \V^\lambda)^2}{d!}\left( \frac{d!}{\dim \V^\lambda} \right)^2
                h_{2g-2}(\lambda),
            \end{equation*}

        \noindent
        expresses $H_g^d$ as the dimension-normalized regular character of a
        uniquely determined central element in the 
        group algebra $\mathbb{C}\group{S}^d$, and we will describe this element explicitly.
        The commutator sum

            \begin{equation*}
                M=\sum_{\pi_1,\pi_2 \in \group{S}^d} \pi_1\pi_2\pi_1^{-1}\pi_2^{-1}
            \end{equation*}

        \noindent
        is a central element in $\C\group{S}^d$ whose eigenvalue in the irreducible representation
        $\V^\lambda$ is $(\frac{d!}{\dim \V^\lambda})^2.$ The Jucys-Murphy elements in $\C\group{S}^d,$

            \begin{equation*}
                J_k = (1\ k)+(2\ k) + \dots + (k-1\ k), \quad 1 \leq k \leq d,
            \end{equation*}

        \noindent
        are not central, but they commute with one another, and for any symmetric polynomial $f$ the element
        $f(J_1,\dots,J_d)$ is central and acts in $\V^\lambda$ as multiplication by $f(\lambda),$ the evaluation 
        of $f$ on the multiset of contents of $\lambda.$
        Thus, $H_g^d$ is the normalized character of $Mh_{2g-2}(J_1,\dots,J_d)$ acting in the regular
        representation. This is also the coefficient of the identity permutation in the class expansion of $Mh_{2g-2}(J_1,\dots,J_d),$ 
        which is the number of configurations $(\pi_1,\pi_2,\tau_1,\dots,\tau_{2g-2})$ of permutations from $\group{S}^d$ satisfying
        the conditions in the theorem statement.
    \end{proof}

    \subsection{Stable limits}
    We now consider the power series

        \begin{equation}
            K(q,t) = 1 + \sum_{d=1}^\infty \frac{q^d}{d!}\sum_{g=1}^\infty t^{2g-2}H_g^d.
        \end{equation}

    \noindent
    Since the coefficient of $\frac{q^d}{d!}$ in $K(q,t)$ is the Maclaurin 
    expansion of the rational function $K^d(t),$ which has radius of convergence $1/(d-1),$
    the bivariate series $K(q,t)$ has zero biradius of convergence and is only valid as a 
    formal power series. This formal power series is the $N=\infty$ limit of the entire function 
    $K_N(q),$ and even though it is not analytic it is useful in understanding 
    the $N \to \infty$ asymptotics of $K_N(q)$ and $L_N(q)=\log K_N(q).$

    Let 

        \begin{equation}
            L(q,t) = \sum_{d=1}^\infty \frac{q^d}{d!}\sum_{g=1}^\infty t^{2g-2}F_g^d
        \end{equation}

    \noindent
    be the unique formal power series such that $K(q,t) = \exp L(q,t)$ in $\C[[q,t]],$
    and let 

        \begin{equation}
            L_N(q) = \sum_{d=1}^\infty \frac{q^d}{d!}L_N^d
        \end{equation}

    \noindent
    be the Maclaurin expansion of the prinicipal logarithm of $K_N(q),$ 
    which has positive radius of convergence. As a consequence of Theorem \ref{thm:LeadingDerivatives}
    we have the following.

        \begin{cor}
            \label{cor:LogarithmicLeadingDerivatives}
            For any $1 \leq d \leq N,$ we have

                \begin{equation*}
                    L_N^d = \sum_{g=1}^\infty N^{2-2g}F_g^d,
                \end{equation*}

            \noindent
            where the series is absolutely convergent. 
        \end{cor}

    Corollary \ref{cor:LogarithmicLeadingDerivatives} says is that if we can 
    compute the formal logarithm of $L(q,t) = \log K(q,t),$ then we will have also computed the 
    stable derivatives $L_N^1,\dots,L_N^N$ of the analytic logarithm $L_N(q)=\log K_N(q).$
    Note that this statement is valid more generally: if 
    $\tilde{K}(q,t)$ is obtained from $K(q,t)$ by formal power series operations,
    and $\tilde{K}_N(q)$ is obtained from $K_N(q)$ by the same sequence of operations,
    then the stable derivatives $\tilde{K}_N^1,\dots,\tilde{K}_N^N$ of 
    $\tilde{K}_N(q)$ are given by the terms of $q$-degree at most $N$ in $\tilde{K}(q,t).$

        \begin{thm}[\cite{GGN1,GGN2}]
            \label{thm:MonotoneConnected}
            The coefficients $F_g^d$ count configurations 
            $(\pi_1,\pi_2,\tau_1,\dots,\tau_{2g-2})$ of permutations from 
            $\group{S}^d$ which satisfy the conditions of Theorem \ref{thm:Monotone},
            and furthermore have the property that $\pi_1,\pi_2,\tau_1,\dots,\tau_{2g-2}$ generate 
            a transitive subgroup of $\group{S}^d.$
        \end{thm}

    An important point is that Theorem \ref{thm:MonotoneConnected} does
    not follow from the Exponential Formula:
    the generating function $K(q,t)$ is exponential in $q$
    but ordinary in $t.$ The fact that the connected/disconnected relationship between $L(q,t)$ and 
    $K(q,t)$ does indeed hold is a consequence of Condition (3) in Theorem \ref{thm:Monotone},
    which prohibits arbitrary shuffling of the transposition 
    factors $\tau_1,\dots,\tau_{2g-2}.$ The fact that $L(q,t)$ is the connected version 
    of $K(q,t)$ even though $t$ is an ordinary rather than exponential marker for the number of transpositions 
    is crucial, since $t$ corresponds to $1/N$ in the asymptotic expansion of unitary 
    matrix integrals, as in Theorem \ref{thm:LeadingDerivatives} above and \cite{MatNov,Novak:BCP}
    more generally. The perturbative parameter $1/N$ does not come with a symmetry
    factor, and the monotonicity condition on the transpositions $\tau_1,\dots,\tau_{2g-2}$ compensates for this.

    \subsection{Hurwitz theory}
    The numbers $H_g^d$ and $F_g^d$ which describe the $N \to \infty$ asymptotics of the $d$th derivatives
    of $K_N(q)$ and $L_N(q)=\log K_N(q)$ at $q=0$ have a topological interpretation which we now explain.
    Let $\mathrm{Conf}_g^d$ be the set of configurations $(\pi_1,\pi_2,\tau_1,\dots,\tau_{2g-2})$
    of permutations from $\group{S}^d$ such that $\pi_1\pi_2\pi_1^{-1}\pi_2^{-1}\tau_1\dots\tau_{2g-2}$
    equals the identity permutation and $\tau_1,\dots,\tau_{2g-2}$ are transpositions. Let 
    $\mathrm{ConConf}_g^d$ be the subset of $\mathrm{Conf}_g^d$ consisting of 
    those configurations such that $\pi_1,\pi_2,\tau_1,\dots,\tau_{2g-2}$
    generate a transitive subgroup of $\group{S}^d.$ Write $B_g^d = |\mathrm{Conf}_g^d|$ and
    $C_g^d=|\mathrm{ConConf}_g^d|.$ 
    By the same argument as in Theorem \ref{thm:Monotone},
    we have 

        \begin{equation}
        \label{eqn:ClassicalContents}
            B_g^d = d!\sum_{\lambda \in \Y^d} \left(\sum_{\Box \in \lambda} c(\Box) \right)^{2g-2},
        \end{equation}

    \noindent
   and by the standard Exponential Formula the generating functions

        \begin{equation}
            B(q,t) = 1 + \sum_{d=1}^\infty \frac{q^d}{d!}\sum_{g=1}^\infty 
            \frac{t^{2g-2}}{(2g-2)!}B_g^d
            \quad\text{and}\quad
            C(q,t) = \sum_{d=1}^\infty \frac{q^d}{d!}\sum_{g=1}^\infty 
            \frac{t^{2g-2}}{(2g-2)!}C_g^d
        \end{equation}

    \noindent
    satisfy 
    
        \begin{equation}
        \label{eqn:ClassicalExponential}
        B(q,t)=\exp C(q,t).
        \end{equation}

    \noindent
    Furthermore, it is a theorem of
    Dijkgraaf \cite{Dijkgraaf} and Kaneko-Zagier \cite{KanZag} that 
    $C_g(q)=\sum_{d=1}^\infty \frac{q^d}{d!} C_g^d$ is a quasimodular form for $g \geq 2.$ The most efficient way to 
    establish this combines \eqref{eqn:ClassicalContents} and \eqref{eqn:ClassicalExponential}
    with the Bloch-Okounkov theorem on quasimodularity
    of arbitrary $q$-uniform averages of shifted symmetric functions \cite{BlOk,CheMolZag,Zagier}.
    
    The symmetric group $\group{S}^d$ acts on both $\mathrm{Conf}_g^d$
    and $\mathrm{ConConf}_g^d$ by simultaneous conjugation, and $\mathrm{Conf}_g^d/\group{S}^d$ 
    and $\mathrm{ConConf}_g^d/\group{S}^d$ are in bijection with
    simple degree $d$ branched covers of an elliptic curve by a compact Riemann surface of
    genus $g.$ More precisely, $|\mathrm{Conf}_g^d/\group{S}^d|$ counts isomorphism classes
    of covers by all genus $g$ surfaces including the disconnected ones, while $|\mathrm{ConConf}_g^d/\group{S}^d|$
    counts connected covers; the genus of a disconnected surface is defined by 
    additivity of the Euler characteristic. This correspondence originates in work of Hurwitz and 
    can be formulated for arbitrary base curves. See \cite{Roth} for elliptic base curves, 
    \cite{EEHS} for the Riemann sphere, and \cite{CavMil} for a general treatment. 
    
    Unfortunately, we do not know how to compute
    $|\mathrm{Conf}_g^d/\group{S}^d|$ and $|\mathrm{ConConf}_g^d/\group{S}^d|.$
    What we do know how to compute is

        \begin{equation}
            \frac{1}{d!}B_g^d = \sum_{\mathcal{O} \in \mathrm{Conf}_g^d} \frac{1}{|\Aut \mathcal{O}|}
            \quad\text{and}\quad
            \frac{1}{d!}C_g^d = \sum_{\mathcal{O} \in \mathrm{ConConf}_g^d} \frac{1}{|\Aut \mathcal{O}|},
        \end{equation}

    \noindent
    where $|\Aut \mathcal{O}|$ is the order of the stabilizer of any particular 
    point in $\mathcal{O}.$ These weighted sums are Hurwitz numbers, and 
    they are generally accepted as computable and useful 
    replacements for the geometrically more natural quantities $|\mathrm{Conf}_g^d/\group{S}^d|$ and $|\mathrm{ConConf}_g^d/\group{S}^d|.$
    Monotone Hurwitz numbers, defined in \cite{GGN1,GGN2} for the sphere and \cite{HIL}
    for arbitrary base curves, ask one to be even more liberal in terms of what is considered acceptable 
    cover counting. For any $\mathcal{O} \in \mathrm{Conf}_g^d/\group{S}^d$ or 
    $\mathrm{ConConf}_g^d/\group{S}^d$, let $\mathrm{Mon}\, \mathcal{O}$
    denote the set of configurations $(\pi_1,\pi_2,\tau_1,\dots,\tau_{2g-2})$ in $\mathcal{O}$ 
    such that the transpositions $\tau_1,\dots,\tau_{2g-2}$ satisfy the monotonicity condition 
    in Theorems \ref{thm:Monotone} and \ref{thm:MonotoneConnected}.
    We then have

        \begin{equation}
            \frac{1}{d!}H_g^d = \sum_{\mathcal{O} \in \mathrm{Conf}_g^d}\frac{1}{|\Aut \mathcal{O}|} \frac{|\mathrm{Mon}\, \mathcal{O}|}{|\mathcal{O}|}
            \quad\text{and}\quad
            \frac{1}{d!}F_g^d = \sum_{\mathcal{O} \in \mathrm{ConConf}_g^d}\frac{1}{|\Aut \mathcal{O}|} 
            \frac{|\mathrm{Mon}\, \mathcal{O}|}{|\mathcal{O}|}.
        \end{equation}

    \noindent
    These sums are disconnected and connected monotone Hurwitz numbers with torus target, 
    and they count the same class of covers as $\frac{1}{d!}B_g^d$ and $\frac{1}{d!}C_g^d$ in an even more non-uniform way.
    Accepting this topologically unnatural cover counting is a reasonable decision insofar as it is both
    computable and useful: every explicit formula for classical Hurwitz numbers has an equally explicit monotone counterpart \cite{GGN1,GGN2}, 
    and monotone Hurwitz numbers occur naturally in asymptotic analysis and mathematical physics. 
    In fact, monotone cover counting is an improvement over previous schemes attempting to compute asymptotic
    expansions of unitary matrix integrals or Yang-Mills partition functions topologically  \cite{Collins,GrossTayl},
    which rely on signed counts of branched covers with contrived branch locus and ramification 
    data, making them even farther from traditional cover counting and much less analytically advantageous.

    \subsection{Concentration inequality}
    For each $g \in \N,$ let 

        \begin{equation}
            F_g(q) = \sum_{d=1}^\infty \frac{q^d}{d!}F_g^d
        \end{equation}

    \noindent
    be the generating function for connected genus $g$ monotone Hurwitz numbers with 
    torus target.
    As explained in Section \ref{sec:FirstOrder}, the series $F_1(q)$ is the logarithm of
    Euler's partition generating function, and quasimodularity of $F_g(q)$ for 
    $g \geq 2$ was established in \cite{HIL}. Thus, for each $m \in \N$ the 
    topologically normalized integral

    \begin{equation}
    \label{eqn:Concentration}
        K_{mN}(q) = e^{-\sum_{g=1}^m F_g(q)}\iint e^{qN\Tr AUA^*U^*} \mathrm{d}U\mathrm{d}A
    \end{equation}

    \noindent
    is an analytic function on $\D.$ We now establish an 
    asymptotic concentration inequality for this normalized integral.

    \begin{thm}
    \label{thm:Concentration}
        For any $\delta \in [0,e^{-1})$ and $m \in \N$ we have 

            \begin{equation*}
                \|K_{mN}-1\|_\delta \leq C_m(\delta) N^{2-2m}
            \end{equation*}

        \noindent
        for all $N \geq N(\delta)$ sufficiently large, where $C_m(\delta) >0$ has no dependence on $N$.
    \end{thm}

    \begin{proof}
    Fix $m \in \N$ and $0 \leq \delta < \gamma < e^{-1}.$
    Consider first the stable $N=\infty$ version of the function $K_{mN}(q)$, i.e. 
    the formal power series

        \begin{equation*}
            K_m(q,t) = e^{-\sum_{g=1}^m t^{2g-2}F_g(q)}K(q,t).
        \end{equation*}

    \noindent
    By Theorem \ref{thm:MonotoneConnected}, this is 

        \begin{equation*}
            K_m(q,t) = e^{-\sum_{g=1}^m t^{2g-2}F_g(q)} e^{\sum_{g=1}^\infty t^{2g-2}F_g(q)} = e^{\sum_{g=m+1}^\infty t^{2g-2}F_g(q)},
        \end{equation*}

    \noindent
    so that 

        \begin{equation*}
            K_m(q,t) = 1 + \sum_{d=1}^\infty \frac{q^d}{d!}\sum_{g=m+1}^\infty t^{2g-2} H_{g \underline{m+1}}^d
        \end{equation*}

    \noindent
    where $H_{g \underline{m+1}}^d$ is the monotone count of possibly disconnected covers of genus $g$ all of 
    whose connected components have genus at least $m+1$.

    Now we utilize the above formal computation to bound the modulus of 
    the analytic function $K_{mN}(q)-1$ on the closure of $\D(\delta).$ 
    Let 

            \begin{equation*}
                K_{mN}(q) = 1 + \sum_{d=1}^\infty \frac{q^d}{d!}K_{mN}^d
            \end{equation*}

        \noindent
        be the Maclaurin expansion of $K_{mN}(q)$ on $\D$. In the stable range, $1 \leq d \leq N$, 
        we have

            \begin{equation*}
                K_{mN}^d = \sum_{g=m+1}^\infty N^{2-2g}H_{g\underline{m+1}}^d \leq \sum_{g=m+1}^\infty N^{2-2g}H_g^d,
            \end{equation*}

        \noindent
        where the upper bound is a convergent series. Now, by Theorem \ref{thm:StirlingBound} and 
        Proposition \ref{prop:StirlingBoung} we have

            \begin{equation*}
                H_g^d \leq d!|\Y^d|{d-1+2g-2 \brace d-1} < d!|\Y^d|e^{d-1}(d-1)^{2g-2},
            \end{equation*}

        \noindent
        so that 

            \begin{equation*}
                \sum_{d=1}^N \frac{\delta^d}{d!}K_{mN}^d < \sum_{d=1}^N (\delta e)^d |\Y^d| \sum_{g=m+1}^\infty 
                \left(\frac{d-1}{N}\right)^{2g-2}.
            \end{equation*}

        \noindent
        For each $1 \leq d \leq N$ we have

            \begin{equation*}
                \sum_{g=m+1}^\infty \left(\frac{d-1}{N}\right)^{2g-2}
                \leq (d-1)^{2m}N^{-2m}\sum_{k=0}^\infty \left(\frac{N-1}{N}\right)^{2k} < (d-1)^{2m}N^{1-2m},
            \end{equation*}

        \noindent
        so

            \begin{equation*}
                \sum_{d=1}^N \frac{\delta^d}{d!}K_{mN}^d < N^{1-2m}\sum_{d=1}^N (\delta e)^d |\Y^d|(d-1)^{2m}.
            \end{equation*}

        \noindent
        Thus, the stable part

            \begin{equation*}
                K_{mNN}(q) = 1 + \sum_{d=1}^N \frac{q^d}{d!} K_{mN}^d 
            \end{equation*}

       \noindent
       of $K_{mN}(q)$ satisfies 

            \begin{equation*}
                \left\| K_{mNN}-1\right\|_\delta < C_m(\delta)N^{2-2m},
            \end{equation*}

        \noindent
       where

            \begin{equation*}
                C_m(\delta) = \sup_{d \in \N}\ (\delta e)^d |\Y^d| (d-1)^{2m}
            \end{equation*}

        \noindent
        is finite because $\delta e <1$ and $|\Y^d| (d-1)^{2m}$ has 
        subexponential growth in $d.$

        Now consider the unstable part of $K_{mN}$, i.e. the tail sum

            \begin{equation}
                K_{mN}(q)-K_{mNN}(q) = \sum_{d=N+1}^\infty \frac{q^d}{d!} K_{mN}^d. 
            \end{equation}

        \noindent
        We have

            \begin{equation*}
                \|K_{mN}\|_\gamma \leq e^{\sum_{g=1}^m \|F_g\|_\gamma}\|K_N\|_\gamma,
            \end{equation*}

        \noindent
        and we have already shown (Proposition \ref{prop:Uniform}) that $\sup_{N\in \N}\|K_N\|_\gamma$
        is finite. Hence 

            \begin{equation*}
                R_m(\gamma) := \sup_{N \in \N} \|K_{mN}\|_\gamma < \infty.
            \end{equation*}

        \noindent
        Thus by Cauchy's inequality we have

            \begin{equation}
                \|K_{mN}(q)-K_{mNN}(q)\|_\delta \leq R_m(\gamma)\sum_{d=N+1}^\infty \left(\frac{\delta}{\gamma} \right)^d 
                =\left(\frac{\delta}{\gamma}\right)^N \frac{R_m(\gamma)\frac{\delta}{\gamma}}{1-\frac{\delta}{\gamma}},
            \end{equation}

        \noindent
        which is exponentially small in $N$.
            
    \end{proof}

\section{Higher Order Asymptotics}
\label{sec:HigherOrder}
In Section \ref{sec:FirstOrder} we proved Theorem \ref{thm:Main} for $m=1$, and 
now we prove the general case using the same 
strategy of deducing convergence of analytic functions
from convergence of their derivatives together with uniform boundedness. 
Fix $m \in \N$ and $0 \leq \varepsilon < \gamma < e^{-1}$.

    \subsection{Convergence of derivatives}
    From the results of Section \ref{sec:FirstOrder},
    there exists $N(\gamma) \in \N$ such that for all $N \geq N(\gamma)$
    the function

        \begin{equation}
            \Delta_{mN}(q) = L_N(q)-\sum_{g=1}^m N^{2-2g} F_g(q)
        \end{equation}

    \noindent
    is defined and analytic in an open neighborhood of the closure of $\D(\gamma)$.
    Henceforth we always assume $N \geq N(\gamma)$, i.e. view $\Delta_{mN}(q)$ as 
    a sequence of analytic functions indexed by $N \geq N(\gamma),$ each of which is 
    defined on an open neighborhood of the closure of $\D(\gamma).$ In particular,
    the Maclaurin expansion

        \begin{equation}
            \Delta_{mN}(q) = 1 + \sum_{d=1}^\infty \frac{q^d}{d!}\Delta_{mN}^d
        \end{equation}

    \noindent
    converges uniformly absolutely on the closure of $\D(\gamma).$ 

        \begin{thm}
            For each $d \in \N$, we have
            
                \begin{equation*}
                        \lim_{N \to \infty} N^{2m-2}\Delta_{mN}^d=0.
                \end{equation*}
        \end{thm}

        \begin{proof}
        By Corollary \ref{cor:LogarithmicLeadingDerivatives},
        for all $N \geq d$ we have

                \begin{equation*}
                    \Delta_{mN}^d = \sum_{g=1}^\infty N^{2-2g}F_g^d - \sum_{g=1}^m N^{2-2g}F_g^d = 
                    \sum_{g=m+1}^\infty N^{2-2g}F_g^d,
                \end{equation*}

            \noindent
            where the series is convergent.
            Thus, for each $N \geq d$ we have $0 \leq N^{2m-2}|\Delta_{mN}^d| \leq N^{-2}R_m^d$, where 

                \begin{equation*}
                    R_m^d = F_{m+1}^d + d^{-2}F_{m+2}^d + d^{-4}F_{m+3}^d + \dots
                \end{equation*}

            \noindent
            is finite and independent of $N$, and we conclude $\lim_{N \to \infty} N^{2m-2}\Delta_{mN}^d=0.$
        \end{proof}

    \subsection{Uniform boundedness}
    We now establish uniform boundedness of 
    $N^{2m-2}\Delta_{mN}(q)$.

        \begin{thm}
            We have

                \begin{equation*}
                    \sup\limits_{N \geq N(\gamma)} N^{2m-2}\|\Delta_{mN}\|_\varepsilon < \infty.
                \end{equation*}
        \end{thm}

        \begin{proof}
            For every $N \geq N(\gamma),$ we have

                \begin{equation*}
                    \|e^{\Delta_{mN}}-1\|_\varepsilon \leq \|e^{\Delta_{mN}}-1\|_\gamma = \|K_{mN}-1\|_\gamma \leq C_m(\gamma)N^{2-2m},
                \end{equation*}

            \noindent
            by Theorem \ref{thm:Concentration}.
            Thus, for any $q \in \C$ with $|q| \leq \gamma$ we have 

                \begin{equation*}
                   \left| e^{\Delta_{mN}(q)} \right| = e^{\Re \Delta_{mN}(q)} \leq 1+C_m(\gamma)N^{2-2m},
                \end{equation*}

            \noindent
            so that 

                \begin{equation*}
                    \sup_{|q| \leq \gamma} \Re \Delta_{mN}(q) \leq C_m(\gamma)N^{2-2m}.
                \end{equation*}

            \noindent
            By the Borel-Carath\'eodory inequality \cite[Section 5.5]{Titchmarsh},

                \begin{equation*}
                    \|\Delta_{mN}\|_\varepsilon \leq \frac{2\varepsilon}{\gamma-\varepsilon}\sup_{|q| \leq \gamma} \Re \Delta_{mN}(q),
                \end{equation*}

            \noindent
            so we obtain 

                \begin{equation*}
                    \|\Delta_{mN}\|_\varepsilon \leq \frac{2\varepsilon C_m(\gamma)}{\gamma-\varepsilon}N^{2-2m}.
                \end{equation*}
            
        \end{proof}

    \subsection{Uniform convergence}
    We now conclude the proof of Theorem \ref{thm:Main}.

        \begin{thm}
            We have $\lim_{N \to \infty} N^{2m-2}\|\Delta_{mN}\|_\varepsilon = 0$.
        \end{thm}

        \begin{proof}
            The proof is identical to the proof of Theorem \ref{thm:ConvergenceToEuler}, with 
            $N^{2m-2}\Delta_{mN}$ replacing $K_N-H_1.$
        \end{proof}

\end{document}